\documentclass[11pt,a4paper]{amsart}

\usepackage{mymacros_english_paper}
\usepackage{mymacros_common}

\title[Nygaard grading, Frobenius--Witt, and regularity]{The Nygaard graded piece, the Frobenius--Witt cotangent complex and regularity criteria}

\author[R. Ishizuka]{Ryo Ishizuka}
\address{Department of Mathematics, Institute of Science Tokyo, 2-12-1 Ookayama, Meguro, Tokyo 152-8551}
\email{ishizuka.r.ac@m.titech.ac.jp}

\author[K. Shimada]{Kanau Shimada}
\address{Department of Mathematics, Institute of Science Tokyo, 2-12-1 Ookayama, Meguro, Tokyo 152-8551}
\email{shimada.k.aj@m.titech.ac.jp}


\begin{document}

\begin{abstract}
For a \(p\)-torsionfree \(p\)-quasisyntomic ring, we construct a fiber sequence relating the \(p\)-th Nygaard graded piece of absolute prismatic cohomology to the Frobenius--Witt cotangent complex. The latter is an arithmetic analogue of the cotangent complex in which the differential of $p$ need not vanish. We also clarify the relationship between the regularity criterion proved by Mao using the Frobenius--Witt cotangent complex and the regularity criterion established by Bhatt and Mathew using the Nygaard graded piece. This relationship was anticipated by Bhatt and Mathew.
\end{abstract}

\maketitle 

\setcounter{tocdepth}{1}
\tableofcontents
\section{Introduction}
Frobenius--Witt derivations are arithmetic analogues of ordinary derivations in which the differential of $p$ need not vanish.
This theory was introduced by T. Saito in \cite{saito2022Frobenius} and by T. Dupuy, E. Katz, J. Rabinoff and D. Zureick-Brown in \cite{dupuy2019Total} (under the name of total $p$-derivations). Furthermore, relationships between such derivations and the regularity of rings were independently clarified by Saito (\cite[Theorem 3.4]{saito2022Frobenius}) and by Hochster and Jeffries (\cite[Corollary 5.10]{HochsterJefferies24}). As one can define the cotangent complex using the derived version of usual derivations, we can also define the Frobenius--Witt cotangent complex using the derived version of Frobenius--Witt derivations (\cite{shimada2026FrobeniusWitt}). We denote this complex by $\FWCot_R\in \mcalD(R/^Lp)$. Using this complex, Z. Mao proved the following regularity criterion:
\footnote{The regularity criterion in \cite[Theorem 5.2]{mao2026FrobeniusWitt} is of a different type from the regularity criteria established in \cite{saito2022Frobenius} and \cite{HochsterJefferies24}. In the latter papers, the regularity criteria are formulated in terms of the freeness and rank of the module of Frobenius--Witt differentials \(\FOmega_R^1\). The fact that \cite[Theorem 5.2]{mao2026FrobeniusWitt} requires only flatness can be regarded as a reflection of the fact that the higher homotopies are taken into account.} 

\begin{theorem}[{\cite[Theorem 5.2]{mao2026FrobeniusWitt}}]\label{theorem1.1}
        Let $R$ be a noetherian ring such that \(p\) belongs to its Jacobson radical $\rad(R)$. Then the following assertions are equivalent:
        \begin{enumerate}
            \item $R$ is regular,
            \item $\FWCot_R$ is flat over $R/^Lp \defeq R \otimes^L_{\setZ} \setF_p$.
        \end{enumerate}
\end{theorem}
On the other hand, B. Bhatt and A. Mathew established the regularity criterion using the theory of prismatic cohomology:
\begin{theorem}[{\cite[Proposition 4.19]{bhatt2023Syntomic}}]\label{theorem1.2}
    Let $R$ be a complete intersection local noetherian ring with residue field $k$ of characteristic $p$. Then the following assertions are equivalent:
    \begin{enumerate}
        \item $R$ is regular,
        \item $\Cofib(\theta_R)\otimes_{R/p}^Lk$ is concentrated in cohomological degree $\ge 1$,
    \end{enumerate}
    where $\theta_R:R/^L p\to (\gr_N^p\prism_R)/^Lp$ is a map of $\mcalD(R/^L p)$ defined in \cite{bhatt2023Syntomic}.
\end{theorem}
In their paper, Bhatt and Mathew expected this theorem to be related to the works of \cite{saito2022Frobenius} and \cite{HochsterJefferies24}. The purpose of this paper is to clarify a relationship between the Frobenius--Witt cotangent complex and the complex $\Cofib(\theta_R)$. More precisely, we prove the following theorem, which connects the Nygaard graded piece of prismatic cohomology and the Frobenius--Witt cotangent complex.
\begin{maintheorem}[\Cref{corollary2.9}]
    Let $R$ be a $p$-torsionfree quasisyntomic ring. Then there is a natural fiber sequence:
    \[\Fib(\alpha(\Cot_{(R/p)/\setF_p})[-p])\to \Cofib(\theta_R)\to \FWCot_R[-1]\to\] in $\mcalD(R/p)$.
\end{maintheorem}
Here, the map $\alpha(\Cot_{(R/p)/\setF_p}) \colon \bigwedge^p \Cot_{(R/p)/\setF_p} \to F_{R/p}^*\Cot_{(R/p)/\setF_p}[p-1]$ is defined in \cite{petrov2025Nondecomposability} by linear algebra.
This map connects the map \(c_R \colon \bigwedge^p \Cot_{(R/p)/\setF_p} \to (R/p)[p]\) defined by using the Sen operator on the diffracted Hodge complex \(\widehat{\Omega}^{\slashed{D}}_R\) and the obstruction map \(\ob_{R, F} \colon F^*_{R/p}\Cot_{(R/p)/\setF_p} \to (R/p)[p]\) for Frobenius liftings constructed by Illusie \cite{illusie1971Complexe}; see \cite{petrov2025Nondecomposability}*{Theorem 7.1} and \Cref{lemma2.8}.

As an application of this theorem, we prove the following.
\begin{maintheorem}[\Cref{theorem2.10}]\label{theorem1.6}
    Let $R$ be a $p$-torsionfree quasisyntomic local ring with residue field $k$. Then the following assertions are equivalent:
    \begin{enumerate}
        \item $\FWCot_R\otimes^L_{R/p}k$ is concentrated in degree zero,
        \item $\Cofib(\theta_R:R/p\to \gr^p_N\prism_R/^Lp)\otimes_{R/p}^Lk$ is concentrated in cohomological degree $\ge 1$.
    \end{enumerate}
\end{maintheorem}
By this theorem, we observe that the regularity criterion via the Frobenius--Witt cotangent complex ($(2)\Rightarrow (1)$ of \Cref{theorem1.1}) formally follows from the regularity criterion via prismatic cohomology ($(2)\Rightarrow (1)$ of \Cref{theorem1.2}) under the assumption that $R$ is a \(p\)-torsionfree complete intersection local noetherian ring with residue characteristic $p$. 

\subsection*{Notation}\label{notation}
Unless otherwise stated, $R$ is assumed to be a $p$-torsionfree quasisyntomic commutative $\setZ_p$-algebra, where $p$ is a fixed prime number (\cite{bhatt2019Topological}*{Definition 4.10}).
We call quasi-categories $\infty$-categories. For any animated ring $A$, we denote the derived $\infty$-category of $A$ by $\mcalD(A)$. We denote by $\widehat{\mcalD}(A)$ the full subcategory of $\mcalD(A)$ consisting of derived $p$-complete objects. The left adjoint of the inclusion $\widehat{\mcalD}(A)\hookrightarrow \mcalD(A)$ is denoted by $M\mapsto \widehat{M}$. We denote the $i$-th derived wedge product of $M\in \widehat{\mcalD}(A)$ by $\bigwedge^iM$ in which we implicitly take the derived $p$-completion.

\subsection*{Acknowledgments}
The authors thank Ryo Suzuki for asking an insightful question about a relationship between \cite{shimada2026Arithmetic}
and \cite{petrov2025Nondecomposability}.

The first-named author was supported by JSPS KAKENHI Grant number 24KJ1085.

\subsection*{AI disclosure}
This work began when the first-named author asked an LLM whether there was a relationship between the Frobenius--Witt cotangent complex and Nygaard filtration of absolute prismatic cohomology.
Through discussions with an LLM, the authors could find the fiber sequence connecting them in one of the main theorems (\Cref{corollary2.9}).
The authors are responsible for verifying the correctness of the statements and proofs in this paper.

\section{Preliminaries}
As is mentioned in Notation, we assume that $R$ is a quasisyntomic $p$-torsionfree ring.
\begin{proposition}[{cf. \cites{bhatt2022Absolute,petrov2025Nondecomposability}}]\label{remark3.2}
    We recall basic properties of the diffracted Hodge complex needed in the sequel.
    \begin{enumerate}
        \item We have a natural equivalence\[\gr_i^{\conj}\widehat\Omega_R^{\slashed D}\simeq \bigwedge^i\widehat{\Cot}_{R/\setZ_p}[-i] \] in $\widehat{\mcalD}(R)$ for any $i\ge 0$. This is \cite[Construction 4.7.1]{bhatt2022Absolute}.
        \item By construction, the Sen operator $\Theta_R$ is compatible with the conjugate filtration on $\Omega_R^{\slashed D}$. Under the identification of $(1)$, the pair $(\gr_i^{\conj}\widehat\Omega_R^{\slashed D},\Theta_R)$ is equivalent to the pair $(\bigwedge^i\widehat{\Cot}_{R/\setZ_p}[-i],-i)$ as objects in the $\infty$-category $\widehat{\mcalD}_\setN(R):=\Fun(B\setN,\widehat{\mcalD}(R))$. This is \cite[Notation 4.7.2]{bhatt2022Absolute}.
        \item By (2), the composition
        \[
        \Fil_p^{\conj}\widehat\Omega_R^{\slashed D}\xrightarrow[]{\Theta_R+p}\Fil_p^{\conj}\widehat\Omega_R^{\slashed D}\xrightarrow[]{\can}\gr_p^{\conj}\widehat\Omega_R^{\slashed D}\xrightarrow[(1)]{\simeq } \bigwedge^p\widehat{\Cot}_{R/\setZ_p}[-p]
        \]
        is naturally equivalent to the composition\[\Fil_p^{\conj}\widehat\Omega_R^{\slashed D}\xrightarrow[]{\can}\gr_p^{\conj}\widehat\Omega_R^{\slashed D}\xrightarrow[(1)]{\simeq } \bigwedge^p\widehat{\Cot}_{R/\setZ_p}[-p]\xrightarrow[]{p-p(=0)} \bigwedge^p\widehat{\Cot}_{R/\setZ_p}[-p],\] which is a zero map. Therefore, the map $\Theta_R+p$ induces the following map:
        \[(\Theta_R+p)_{\mathrm{low}}:\Fil_p^{\conj}\widehat\Omega_R^{\slashed D}\to \Fil_{p-1}^{\conj}\widehat\Omega_R^{\slashed D}.\] 
        Then there is a natural equivalence in $\widehat{\mcalD}(R)$:\[\gr_N^p\prism _R\simeq \Fib((\Theta_R+p)_{\mathrm{low}}:\Fil_p^{\conj}\widehat\Omega_R^{\slashed D}\to \Fil_{p-1}^{\conj}\widehat\Omega_R^{\slashed D}),\] where $\gr_N^p\prism _R$ is the $p$-th graded piece of the Nygaard filtration on the prismatic cohomology $\prism_R$. This is \cite[Remark 5.5.8]{bhatt2022Absolute}.
        \item By (2), we have a natural decomposition:\[\Fil_p^{\conj}\widehat\Omega_R^{\slashed D}\simeq \bigoplus _{i=0}^{p-1}(\Fil_p^{\conj}\widehat\Omega_{R}^{\slashed D})_i,\] where the endomorphism \(\Theta_R + i\) acts topologically nilpotently on each \((\Fil_p^{\conj}\widehat\Omega_{R}^{\slashed D})_i\).
        Moreover, we have a natural equivalence:
        \[(\Fil_p^{\conj}\widehat\Omega_R^{\slashed D})_i\simeq \bigwedge^i\widehat{\Cot}_{R/\setZ_p}[-i]\] for any $1\le i\le p-1$. This is \cite[Section 6.2]{petrov2025Nondecomposability}.
    \end{enumerate}
\end{proposition}
The following lemma is a direct consequence of the facts recalled in \Cref{remark3.2}.
\begin{lemma}\label{lemma2.3}
    There is a natural fiber sequence:
    \[\gr_N^p\prism _R\to (\Fil_p^{\conj}\widehat\Omega_R^{\slashed D})_0\xrightarrow[]{(\Theta_R+p)_{\mathrm{low}}} R\to\] in $\widehat{\mcalD}(R)$.
\end{lemma}
\begin{proof}
    By (3) of \Cref{remark3.2}, we have the following natural fiber sequence:
    \[\gr_N^p\prism _R\to \Fil_p^{\conj}\widehat\Omega_R^{\slashed D}\xrightarrow[]{(\Theta_R+p)_{\mathrm{low}}} \Fil_{p-1}^{\conj}\widehat\Omega_R^{\slashed D}\to \] in $\widehat{\mcalD}(R)$. By the compatibility of the Sen operator with the conjugate filtration, the map \[(\Theta_R+p)_{\mathrm{low}}:\Fil_p^{\conj}\widehat\Omega_R^{\slashed D}\to \Fil_{p-1}^{\conj}\widehat\Omega_R^{\slashed D}\] can be written as the direct sum of the following maps:
    \[(\Fil_p^{\conj}\widehat\Omega_R^{\slashed D})_0\to (\Fil_{p-1}^{\conj}\widehat\Omega_R^{\slashed D})_0\simeq R\]
    \[\bigoplus_{i=1}^{p-1}(\Fil_p^{\conj}\widehat\Omega_R^{\slashed D})_i\to \bigoplus_{i=1}^{p-1}(\Fil_{p-1}^{\conj}\widehat\Omega_R^{\slashed D})_i.\] Since the latter map can be identified with the following map by (4) of \Cref{remark3.2}:\[\bigoplus_{i=1}^{p-1}\bigwedge^i\widehat{\Cot}_{R/\setZ_p}[-i]\to \bigoplus_{i=1}^{p-1}\bigwedge^i\widehat{\Cot}_{R/\setZ_p}[-i]:(a_1,\dots,a_{p-1})\mapsto ((p-1)a_1,(p-2)a_2,\dots,a_{p-1}),\] the latter map is an equivalence. Therefore, the above fiber sequence induces the desired fiber sequence. 
\end{proof}

Next, we recall the definition of the map $c_R$ defined in \cite{petrov2025Nondecomposability}. To do this, we first recall some preliminary facts on the $\infty$-category $\widehat{\mcalD}_\setN(R)$ in \cite[Section 6.1]{petrov2025Nondecomposability}.
\begin{remark}\label{remark2.4new}
    In \cite{petrov2025Nondecomposability}, Petrov considers the $\infty$-category $\widehat{\mcalD}_\setN(R)$ defined to be:
    \[\widehat{\mcalD}_\setN(R):=\Fun(B\setN,\widehat{\mcalD}(R)).\] We recall some facts on its mapping animae.
    \begin{enumerate}
        \item By \citeKero{00JC}, there is a canonical map $\Delta^1/\partial\Delta^1\to B\setN$, which is an inner anodyne map. Thus, the induced functor \[\widehat{\mcalD}_\setN(R)=\Fun(B\setN,\widehat{\mcalD}(R))\to \Fun(\Delta^1/\partial\Delta^1,\widehat{\mcalD}(R))\] is an equivalence of $\infty$-categories. Therefore, as is written in \cite[Section 6.1]{petrov2025Nondecomposability}, the mapping anima can be computed as follows:
        \begin{align*}
        \Hom_{\widehat{\mcalD}_\setN(R)}((M,f),(N,g))&\simeq \Eq(\Hom_{\widehat{\mcalD}(R)}(M,N)\rightrightarrows \Hom_{\widehat{\mcalD}(R)}(M,N))\\
        &\simeq \Fib(\Hom_{\widehat{\mcalD}(R)}(M,N)\to \Hom_{\widehat{\mcalD}(R)}(M,N):h\mapsto h\circ f-g\circ h),
        \end{align*}
        where the transition maps in the equalizer are defined to be $h\mapsto h\circ f$ and $h\mapsto g\circ h$.
        \item By (1), there is the following canonical commutative diagram of animae:
        \[ \begin{tikzcd}[column sep=small]
        \Hom_{\widehat{\mcalD}_\setN(R)}((M,f),(N,g)) \arrow[r, "\simeq "] \arrow[d, "\can"'] & \Fib(\Hom_{\widehat{\mcalD}(R)}(M,N)\to \Hom_{\widehat{\mcalD}(R)}(M,N)) \arrow[d, "\can"] \\ \Hom_{\mcalD_\setN(R/p)}((M/^Lp,f/^Lp),(N/^Lp,g/^Lp)) \arrow[r, "\simeq"'] & \Fib(\Hom_{\mcalD(R/p)}(M/^Lp,N/^Lp)\to \Hom_{\mcalD(R/p)}(M/^Lp,N/^Lp)) \end{tikzcd} \] 
        \end{enumerate}
\end{remark}
Although the following two lemmas on the $\infty$-category $\widehat{\mcalD}_\setN(R)$ are not necessary for the construction of the map $c_R$, we recall them here for later use.
\begin{lemma}\label{remark2.8}
    Let \[(A,0)\xrightarrow[]{i} (B,f)\xrightarrow[]{q} (C,0)\to \] be a fiber sequence in the $\infty$-category $\widehat{\mcalD}_\setN(R)$. Then by the universality of $A$ as a fiber of the map $q:B\to C$, $f$ induces a natural map $f_{\mathrm
    {low}}:B\to A$. Moreover, by the universality of $C$ as a cofiber of the map $i:A\to B$, $f_\mathrm{low}$ further induces a natural map $h:C\to A$. Then there is a natural equivalence:
    \[f_\mathrm{low}\simeq h\circ q.\]
\end{lemma}
\begin{proof}
    This follows from the commutativity of the diagram \((A,0)\xrightarrow[]{i} (B,f)\xrightarrow[]{q} (C,0)\to\).
\end{proof}
\begin{lemma}\label{remark2.9}
    Keep the notation in \Cref{remark2.8}.
    Let $c:C\to A[1]\oplus A$ be the map classifying the connecting map of this fiber sequence ((1) of \Cref{remark2.4new}). Then the composition 
    \[C\xrightarrow[]{c}A[1]\oplus A\xrightarrow[]{\pr}A\] coincides with the map $h:C\to A$. 
\end{lemma}
\begin{proof}
This follows from the paragraph between \cite[Corollary 6.6]{petrov2025Nondecomposability} and \cite[Lemma 6.7]{petrov2025Nondecomposability}.
\end{proof}

\begin{construction}\label{remark2.4}
    We recall the map
    \[c_R:\bigwedge^p {\Cot}_{(R/p)/\setF_p}\to (R/p)[p]\] in $\mcalD(R/p)$ defined in \cite[Notation 6.4]{petrov2025Nondecomposability}. Consider the following fiber sequence in $\widehat{\mcalD}_{\setN}(R)$:
    \[(\Fil_{p-1}^{\conj}\widehat\Omega_R^{\slashed D},\Theta_R)\to (\Fil_p^{\conj}\widehat\Omega_R^{\slashed D},\Theta_R)\to (\bigwedge^p\widehat{\Cot}_{R/\setZ_p}[-p],-p)\to. \] 
    We have the following equivalences in $\widehat{\mcalD}(R)$:
    \begin{align*}
    &\RHom_{\widehat{\mcalD}_\setN(R)}((\bigwedge^p\widehat{\Cot}_{R/\setZ_p}[-p],-p),(\Fil_{p-1}^{\conj}\widehat\Omega_R^{\slashed D},\Theta_R)[1])\\
    &\simeq \bigoplus_{i=0}^{p-1}\RHom_{\widehat{\mcalD}_\setN(R)}((\bigwedge^p\widehat{\Cot}_{R/\setZ_p}[-p],-p),(\bigwedge^i\widehat{\Cot}_{R/\setZ_p}[-i],-i)[1])\\
    &\simeq \bigoplus_{i=0}^{p-1}\Fib(\RHom_{\widehat{\mcalD}(R)}(\bigwedge^p\widehat{\Cot}_{R/\setZ_p}[-p],\bigwedge^i\widehat{\Cot}_{R/\setZ_p}[-i][1])\xrightarrow[]{\times (i-p)} \RHom_{\widehat{\mcalD}(R)}(\bigwedge^p\widehat{\Cot}_{R/\setZ_p}[-p],\bigwedge^i\widehat{\Cot}_{R/\setZ_p}[-i][1]))\\
    &\simeq \Fib(\RHom_{\widehat{\mcalD}(R)}(\bigwedge^p\widehat{\Cot}_{R/\setZ_p}[-p],R[1])\xrightarrow[]{\times (-p)} \RHom_{\widehat{\mcalD}(R)}(\bigwedge^p\widehat{\Cot}_{R/\setZ_p}[-p],R[1]))\\
    &\simeq \RHom_{\widehat{\mcalD}(R)}(\bigwedge^p\widehat{\Cot}_{R/\setZ_p},R/p[p])\\
    &\simeq \RHom_{\mcalD(R/p)}(\bigwedge^p{\Cot}_{(R/p)/\setF_p},R/p[p]),
    \end{align*}
    where in the second equivalence, we used the identification in (1) of \Cref{remark2.4new}.
    Therefore, the connecting map of the above fiber sequence defines a map:
    \[c_R:\bigwedge^p\Cot_{(R/p)/\setF_p}\to (R/p)[p], \] which is the desired one.
\end{construction}

\begin{construction}\label{remark2.5}
    To state our first main result, we recall the map $\theta_R$ defined in \cite[]{bhatt2023Syntomic}. Since (3) of \Cref{remark3.2} implies the equivalence\[\gr_N^p\prism_{\setZ_p}\simeq \Fib(\setZ_p\xrightarrow[]{p}\setZ_p)\simeq \setF_p[-1],\] we have the natural map:
    \[\setF_p\to \setF_p\oplus\setF_p[-1]\simeq \gr_N^p\prism_{\setZ_p}/^Lp\to\gr_N^p\prism_R/^Lp,\] 
    where the first map is the inclusion into the first factor and the second map is induced by the canonical map $\setZ_p\to R$.
    We denote by $\theta_R$ the $R/p$-linearization
    \[R/p\to \gr_N^p\prism_R/^Lp\] of the above $\setF_p$-linear map.
\end{construction}
At the end of this section, we recall the following lemma, which is a stable $\infty$-categorical version of the octahedral axiom for triangulated categories.
\begin{lemma}\label{lemma2.9new}
    Let $\mcalD$ be a stable $\infty$-category. Let $f:A\to B$ and $g:B\to C$ be maps in $\mcalD$. Then there is a natural fiber sequence:
    \[\Fib(f)\to \Fib(g\circ f)\to \Fib(g)\to \] in $\mcalD$.
\end{lemma}
\begin{proof}
    This is (TR4) in the proof of \cite[Theorem 1.1.2.14]{lurie2017Higher}.
\end{proof}

\section{The existence of the fiber sequence}
As mentioned in Notation, we assume that $R$ is a quasisyntomic $p$-torsionfree ring.
We first identify the cofiber of the map $\theta_R$ with the fiber of the map $c_R[-p]$ (\Cref{proposition2.6}). We then factorize \(c_R\) by the map \(\alpha(\Cot_{(R/p)/\setF_p})\) defined in \cite{petrov2025Nondecomposability} in \Cref{construction2.10new} and the arithmetic Kodaira--Spencer map \(\kappa_R\) defined in \cite{shimada2026Arithmetic} (\Cref{lemma2.8}) and conclude the proof of \Cref{theorem1.1}.

\begin{theorem}\label{proposition2.6}
    There is a natural fiber sequence:
    \[\Cofib(\theta_R)\to \bigwedge^p\Cot_{(R/p)/\setF_p}[-p]\xrightarrow[]{c_R[-p]}R/p\to \] in $\mcalD(R/p)$.
\end{theorem}
\begin{proof}
    Consider the fiber sequence:
    \[(R,p)\to ((\Fil_p^{\conj}\widehat\Omega_R^{\slashed D})_0,\Theta_R+p)\xrightarrow[]{q} (\bigwedge ^p\widehat{\Cot}_{R/\setZ_p}[-p],0)\to \] in $\widehat{\mcalD}_\setN(R)$. Let
    \[c':\bigwedge ^p\widehat{\Cot}_{R/\setZ_p}[-p]\to R/p\] be the map classifying the connecting map of the above fiber sequence (see (1) of \Cref{remark2.4new}). By (2) of \Cref{remark2.4new}, the map
    \[c'/^Lp:\bigwedge ^p\Cot_{(R/p)/\setF_p}[-p]\to R/p\oplus R/p[1]\] obtained by applying $(-)\otimes_R^LR/p$ to the map $c'$ classifies the fiber sequence:
    \begin{equation} \label{eq1}
        (R/p,0)\to ((\Fil_p^{\conj}\widehat\Omega_R^{\slashed D})_0/^Lp,\Theta_R+p)\xrightarrow[]{q} (\bigwedge ^p\Cot_{(R/p)/\setF_p}[-p],0)\to
    \end{equation}
    in $\mcalD_\setN(R/p)$.
    On the other hand, let $c'_R$ be the map classifying the connecting map of the following fiber sequence (see (1) of \Cref{remark2.4new}):
    \[(R,0)\to ((\Fil_p^{\conj}\widehat\Omega_R^{\slashed D})_0,\Theta_R)\to (\bigwedge ^p\widehat{\Cot}_{R/\setZ_p}[-p],-p)\to \] in $\widehat{\mcalD}_\setN(R)$. By definition, the map $c_R[-p]$ in \Cref{remark2.4} is the $R/p$-linearization of $c'_R$. By (2) of \Cref{remark2.4new}, the map $c'_R/^Lp$ classifies the connecting map of the fiber sequence \eqref{eq1}. Combining them, we see that the map $c_R[-p]$ can be written as the following composition:
    \[\bigwedge ^p\Cot_{(R/p)/\setF_p}[-p]\xrightarrow[]{c'/^Lp}R/p\oplus R/p[1]\xrightarrow[]{\pr}R/p.\] By \Cref{remark2.9} and \Cref{remark2.8}, we have the following canonical factorization:
    \begin{equation} \label{eq2}
        (\Theta_R+p)_{\mathrm{low}}/^Lp:(\Fil_p^{\conj}\widehat\Omega_R^{\slashed D})_0/^Lp\xrightarrow[]{q}\bigwedge ^p{\Cot}_{(R/p)/\setF_p}[-p]\xrightarrow[]{c_R[-p]}R/p.
    \end{equation}
    Applying \Cref{lemma2.9new} to this factorization, we get the following natural fiber sequence:
    \[\Fib(q)\to \Fib((\Theta_R+p)_{\mathrm{low}}/^Lp)\to \Fib(c_R[-p])\to. \] By \Cref{lemma2.3}, this fiber sequence can be naturally identified with the following one \[R/p\to \gr_N^p\prism_R/^Lp\to \Fib(c_R[-p])\to .\]
    To complete the proof, it remains to prove that the map $R/p\to \gr_N^p\prism_R/^Lp$ in the above fiber sequence is naturally equivalent to the map $\theta_R$. By the definition of $\theta_R$ (see \Cref{remark2.5}) and the functoriality of the above fiber sequence, we may assume $R=\setZ_p$. In this case, the factorization \eqref{eq2} can be identified with the following composition:
    \[\setF_p\to 0\to \setF_p.\] The first map in the fiber sequence obtained from this composition by applying \Cref{lemma2.9new} can be written as follows:
    \[\setF_p\xrightarrow[]{\can}\setF_p\oplus\setF_p[-1].\] By \Cref{remark2.5}, this coincides with the map $\theta_{\setZ_p}$.
\end{proof}
\begin{remark}
    In the proof of \cite[Proposition 4.19]{bhatt2023Syntomic}, Bhatt and Mathew also constructed the fiber sequence in the form of \Cref{proposition2.6} by animation. This fiber sequence coincides with the one in \Cref{proposition2.6}. More precisely, there exists an automorphism $u_R$ on $\bigwedge^p\Cot_{(R/p)/\setF_p}[-p]$ and a commutative diagram:
    \begin{equation} \label{commudiagram_theta}
        \begin{tikzcd}[column sep=2cm, row sep=1cm]
        \Cofib(\theta_R) \arrow[r] \arrow[d, "\id"] & \bigwedge^p\Cot_{(R/p)/\setF_p}[-p] \arrow[r,"{c_R[-p]}"] \arrow[d,"u_R"] & R/p \arrow[r] \arrow[d,"\cong"] & {} \\
        \Cofib(\theta_R) \arrow[r] & \bigwedge^p\Cot_{(R/p)/\setF_p}[-p] \arrow[r,""] & R/p \arrow[r] & {}
        \end{tikzcd}
    \end{equation} 
    in $\mcalD(R/p)$,
where the upper fiber sequence is the one in \Cref{proposition2.6} and the lower one is the fiber sequence (17) in the proof of \cite[Proposition 4.19]{bhatt2023Syntomic}.

To prove this, we may assume $R$ is a \(p\)-completed polynomial $\setZ_p$-algebra. In this case, the proof of \cite[Proposition 4.19]{bhatt2023Syntomic} implies that $\Cofib(\theta_R)$ is concentrated in cohomological degrees $1$ and $p$ with $H^1\cong R/p$ and $H^p\cong \Omega^p_{(R/p)/\setF_p}$. Therefore, we have a natural automorphism 
\[u_R:\Omega^p_{(R/p)/\setF_p}[-p]\to \Omega^p_{(R/p)/\setF_p}[-p]\] fitting into the following commutative diagram:
\[ \begin{tikzcd} \Cofib(\theta_R) \arrow[r, ""] \arrow[d, "\id"'] & \Omega^p_{(R/p)/\setF_p}[-p] \arrow[d, "u_R"] \\ \Cofib(\theta_R)  \arrow[r, ""'] & \Omega^p_{(R/p)/\setF_p}[-p], \end{tikzcd} \] where the upper horizontal map is the one in the fiber sequence of \Cref{proposition2.6} and the lower horizontal map is the one in the fiber sequence (17) in the proof of \cite[Proposition 4.19]{bhatt2023Syntomic}. This square extends to the desired commutative diagram \eqref{commudiagram_theta}.
\end{remark}

Next, we connect the above result with the Frobenius--Witt cotangent complex of $R$. To do so, we establish the following result (\Cref{lemma2.8}) using the result of \cite{petrov2025Nondecomposability}. For convenience, we recall the construction of the map $\alpha(\Cot_{(R/p)/\setF_p})$ defined in \cite{petrov2025Nondecomposability} and the map $\kappa_R$ defined in \cite{shimada2026Arithmetic}.
The latter map is inspired by the arithmetic Kodaira--Spencer map defined in \cite{dupuy2019Total} in the case of smooth algebras.

\begin{construction}\label{construction2.10new}
    Let $E$ be an object of $\mcalD(R/p)$. We define the map $\alpha(E)$ in $\mcalD(R/p)$ to be the following composition:
    \[\alpha(E):\Lambda^pE\simeq \Gamma^p(E[-1])[p]\xrightarrow[]{\psi_{E[-1]}[p]}F^*_{R/p}E[p-1],\] where the first equivalence is Illusie's d\'{e}calage equivalence (\cite{illusie1971Complexe}*{Proposition I.4.3.2.1}; see also \cite{petrov2025Nondecomposability}*{Lemma 2.6}) and $\psi_{E[-1]}$ is the map defined in the paragraph between \cite[Lemma 2.5]{petrov2025Nondecomposability} and \cite[Lemma 2.6]{petrov2025Nondecomposability}.
\end{construction}
\begin{remark}
    Precisely, the definition of the map $\alpha(E)$ in \cite{petrov2025Nondecomposability} is different from the one in \Cref{construction2.10new}. However, these definitions are the same by the paragraph between \cite[Definition 3.4]{petrov2025Nondecomposability} and \cite[Remark 3.5]{petrov2025Nondecomposability}.
\end{remark}
\begin{construction}\label{construction2.11new}
    By \cite[Theorem 5.3]{shimada2026FrobeniusWitt}, we have the following natural fiber sequence:
    \[R/p\to \FWCot_R\to F^*_{(R/p)}\Cot_{(R/p)/\setF_p}\to \] in $\mcalD(R/p)$. The arithmetic Kodaira--Spencer map $\kappa_R$ is defined to be the connecting map \[F^*_{R/p}\Cot_{(R/p)/\setF_p}\to R/p[1]\] of the above fiber sequence.
\end{construction}

\begin{proposition}\label{lemma2.8}
    The map 
    \[c_R:\bigwedge^p\Cot_{(R/p)/\setF_p}\to (R/p)[p]\]
    naturally factors as follows:
    \[\bigwedge^p\Cot_{(R/p)/\setF_p}\xrightarrow[]{\alpha(\Cot_{(R/p)/\setF_p})}F^*_{R/p}\Cot_{(R/p)/\setF_p}[p-1]\xrightarrow[]{\kappa_R[p-1]}R/p[p]\] in $\mcalD(R/p)$. 
\end{proposition}
\begin{proof}
    By \cite[Theorem 7.1]{petrov2025Nondecomposability}, the map $c_R$ has a natural factorization as follows:
    \[c_R \colon \bigwedge^p\Cot_{(R/p)/\setF_p}\xrightarrow[]{\alpha(\Cot_{(R/p)/\setF_p})}F^*_{R/p}\Cot_{(R/p)/\setF_p}[p-1]\xrightarrow[]{\ob_{F,R}[p-1]}R/p[p],\] where $\ob_{F,R}$ is the map defined in \cite[Convention 5.5 and Lemma 5.6]{petrov2025Nondecomposability}, which comes from the obstruction class defined in \cite{illusie1971Complexe}*{III.2.2}.
    Therefore, it remains to show that there is a natural equivalence \(\kappa_R\simeq \ob_{F,R}\).
    Let $\lambda_R$ be the map defined in \cite[Definition 5.2]{shimada2026Arithmetic}, which coincides with the map $\gamma_{R_0}-\gamma_{R_0^{(1)}}$ in \cite[Lemma 5.6]{petrov2025Nondecomposability}. By \cite[Theorem 5.3]{shimada2026Arithmetic}, the map $\lambda_R$ is naturally equivalent to the following composition:
    \begin{equation} \label{eq:factorization_cR}
        \lambda_R \colon F^*_{R/p}\Cot_{(R/p)/(\setZ/p^2)}\xrightarrow[]{\can}F^*_{R/p}\Cot_{(R/p)/\setF_p}\xrightarrow[]{\kappa_R}(R/p)[1].
    \end{equation}
    On the other hand, the map $\ob_{F,R}$ is defined as the following composition:
    \[\ob_{R, F} \colon F^*_{R/p}\Cot_{(R/p)/\setF_p}\xrightarrow[]{F^*_{R/p}(s)}F^*_{R/p}\Cot_{(R/p)/(\setZ/p^2)}\xrightarrow[]{\lambda_R}(R/p)[1],\] where $F^*_{R/p}(s)$ is the Frobenius pullback of the following natural composition:
    \begin{equation} \label{eq:factorization_cR2}
        s \colon \Cot_{(R/p)/\setF_p}\simeq \Cot_{(R/p^2)/(\setZ/p^2)}\otimes^L_{R/p^2}R/p\xrightarrow[]{\can}\Cot_{(R/p)/(\setZ/p^2)}.
    \end{equation}
    By definition, the map $F^*_{R/p}(s)$ is a section of the canonical map\[\can:F^*_{R/p}\Cot_{(R/p)/(\setZ/p^2)}\to F^*_{R/p}\Cot_{(R/p)/\setF_p}.\] From this, \eqref{eq:factorization_cR} and \eqref{eq:factorization_cR2}, we get the desired natural identification \(\kappa_R\simeq \ob_{F,R}\).
\end{proof}

\begin{theorem}\label{corollary2.9}
    There is a natural fiber sequence:
    \[\Fib(\alpha(\Cot_{(R/p)/\setF_p})[-p])\to \Cofib(\theta_R)\to \FWCot_R[-1]\to\] in $\mcalD(R/p)$.
\end{theorem}
\begin{proof}
    By \Cref{lemma2.8}, there is a natural factorization:
    \[c_R[-p]\simeq \kappa_R[-1]\circ \alpha(\Cot_{(R/p)/\setF_p})[-p].\]
    Applying \Cref{lemma2.9new} to this factorization, we have the following fiber sequence:
    \[\Fib(\alpha(\Cot_{(R/p)/\setF_p})[-p])\to \Fib(c_R[-p])\to \Fib(\kappa_R[-1])\to.\]
    On the other hand, \Cref{proposition2.6} implies \[\Fib(c_R[-p])\simeq \Cofib(\theta_R).\]
    By the definition of the arithmetic Kodaira--Spencer map, we have 
    \[\Fib(\kappa_R[-1])\simeq \FWCot_R[-1].\]
    Combining them, we get the desired fiber sequence.
\end{proof}
\section{The comparison of regularity criteria}
Using the results in the previous section, we study the relationship between the regularity criterion using the Frobenius--Witt cotangent complex and the one established in \cite{bhatt2023Syntomic}.
\begin{theorem}\label{theorem2.10}
    Let $R$ be a $p$-torsionfree quasisyntomic local ring with residue field $k$. Then the following assertions are equivalent:
    \begin{enumerate}
        \item $\FWCot_R\otimes^L_{R/p}k$ is concentrated in degree zero,
        \item $\Cofib(\theta_R:R/p\to \gr^p_N\prism_R/^Lp)\otimes_{R/p}^Lk$ is concentrated in cohomological degree $\ge 1$.
    \end{enumerate}
\end{theorem}
\begin{proof}
    By definition and \Cref{lemma2.8}, we have the following commutative diagram of fiber sequences in $\mcalD(k)$:
\[\begin{tikzcd}[column sep=2cm, row sep=1cm]
\Cofib(\theta_R)\otimes^L_{R/p}k \arrow[r] \arrow[d] & \bigwedge^p\Cot_{(R/p)/\setF_p}[-p]\otimes^L_{R/p}k \arrow[r,"{c_R[-p]\otimes^L k}"] \arrow[d,"{\alpha(\Cot_{(R/p)/\setF_p})[-p]\otimes ^L k}"'] & k \arrow[r] \arrow[d,"\id"] & {} \\
\FWCot_R[-1]\otimes^L_{R/p}k \arrow[r] & (F^*_{R/p}\Cot_{(R/p)/\setF_p}[-1])\otimes^L_{R/p}k \arrow[r,"{\kappa_R[-1]\otimes^L k}"] & k \arrow[r] & {,}
\end{tikzcd}\]
where the upper fiber sequence is the one in \Cref{proposition2.6} and the lower fiber sequence is obtained by applying $(-)\otimes^L_{R/p}k$ to the fundamental fiber sequence of Frobenius--Witt cotangent complexes (\cite[Theorem 5.3]{shimada2026FrobeniusWitt}). As is mentioned in the proof of \cite[Proposition 4.19]{bhatt2023Syntomic}, the assertion (2) is equivalent to the injectivity of the map \[H_0(c_R[-p]\otimes^L k):H_0(\bigwedge^p\Cot_{(R/p)/\setF_p}[-p]\otimes^L_{R/p}k)\to k.\] Based on this, we first compute $H_0(\bigwedge^p\Cot_{(R/p)/\setF_p}[-p]\otimes^L_{R/p}k)$:
Since $\Cot_{(R/p)/\setF_p}\otimes_{R/p}^Lk$ is concentrated in homological degree $[0,1]$, we have the following fiber sequence:
\[H_1(\Cot_{(R/p)/\setF_p}\otimes_{R/p}^Lk)[1]\to \Cot_{(R/p)/\setF_p}\otimes_{R/p}^Lk\to H_0(\Cot_{(R/p)/\setF_p}\otimes_{R/p}^Lk)\to .\] Thus, we have a finite filtration $\Fil_{\mathrm{wedge}}^\bullet$ on $\bigwedge^p(\Cot_{(R/p)/\setF_p}\otimes_{R/p}^Lk)$ whose graded pieces are \begin{align*}
    \gr^i_{\mathrm{wedge}}&\simeq (\bigwedge^i(H_1(\Cot_{(R/p)/\setF_p}\otimes_{R/p}^Lk)[1]))\otimes_k^L(\bigwedge^{p-i}H_0(\Cot_{(R/p)/\setF_p}\otimes_{R/p}^L k) \\
    &\simeq (\Gamma^i(H_1(\Cot_{(R/p)/\setF_p}\otimes_{R/p}^Lk))[i])\otimes_k^L(\bigwedge^{p-i} H_0(\Cot_{(R/p)/\setF_p}\otimes_{R/p}^Lk)) 
\end{align*}
for $0\le i\le p$, where in the second equivalence, we used Illusie's d\'{e}calage isomorphism (e.g. \cite[Lemma 2.6]{petrov2025Nondecomposability}).
Using this, we see that the natural morphism \(H_1(\Cot_{(R/p)/\setF_p} \otimes^L_{R/p} k)[1] \to \Cot_{(R/p)/\setF_p} \otimes^L_{R/p} k\) induces a natural equivalence
\begin{equation*}
    \Gamma^p(H_1(\Cot_{(R/p)/\setF_p}\otimes^L_{R/p}k)) \xrightarrow{\simeq} H_p(\bigwedge^p(\Cot_{(R/p)/\setF_p}\otimes^L_{R/p}k) \simeq H_0(\bigwedge^p\Cot_{(R/p)/\setF_p}[-p]\otimes^L_{R/p}k),
\end{equation*}
Moreover, the map
\[H_0(\alpha(\Cot_{(R/p)/\setF_p})[-p]\otimes^L k): H_0(\bigwedge^p\Cot_{(R/p)/\setF_p}[-p]\otimes^L_{R/p}k)\to H^0((F^*_{R/p}\Cot_{(R/p)/\setF_p}[-1])\otimes^L_{R/p}k)\] can be identified with the following map:
\[\psi_{H_1(\Cot_{(R/p)/\setF_p}\otimes^L_{R/p}k)}:\Gamma^p(H_1(\Cot_{(R/p)/\setF_p}\otimes^L_{R/p}k))\to F^*_kH_1(\Cot_{(R/p)/\setF_p}\otimes_{R/p}^Lk),\] which is defined in the paragraph between \cite[Lemma 2.5]{petrov2025Nondecomposability} and \cite[Lemma 2.6]{petrov2025Nondecomposability}. Indeed, by the definition of $\alpha(\Cot_{(R/p)/\setF_p})$, we have 
\[\alpha(\Cot_{(R/p)/\setF_p})[-p]\otimes_{R/p}^Lk\simeq \psi_{\Cot_{(R/p)/\setF_p}[-1]}\otimes_{R/p}^Lk.\]
By the definition of $\psi$, we have a natural equivalence:
\[\psi_{\Cot_{(R/p)/\setF_p}[-1]}\otimes^L_{R/p}k\simeq \psi_{\Cot_{(R/p)/\setF_p}\otimes^L_{R/p}k[-1]}.\] Thus, we have 
\[H^0(\alpha(\Cot_{(R/p)/\setF_p})[-p]\otimes_{R/p}^Lk)\simeq H_0(\psi_{\Cot_{(R/p)/\setF_p}\otimes^L_{R/p}k[-1]}).\] By the functoriality of $\psi_{(-)}$, there is a natural morphism \[\psi_{H_1(\Cot_{(R/p)/\setF_p}\otimes^L_{R/p}k)}\to \psi_{\Cot_{(R/p)/\setF_p}\otimes^L_{R/p}k[-1]}.\] Therefore, we get the desired identification:
\[H_0(\alpha(\Cot_{(R/p)/\setF_p})[-p]\otimes^L k)\cong \psi_{H_1(\Cot_{(R/p)/\setF_p}\otimes^L_{R/p}k)}.\]
Using these observations, we will prove the equivalence of (1) and (2).

$(1)\Rightarrow (2)$: Assume that $\FWCot_R\otimes_{R/p}^Lk$ is concentrated in degree zero. Then the map
\[H^0(\kappa_R[-1]\otimes^L k):F^*_kH_1(\Cot_{(R/p)/\setF_p}\otimes_{R/p}^Lk)\to k\] is injective. In particular, we have
\[\dim_kH_1(\Cot_{(R/p)/\setF_p}\otimes_{R/p}^Lk)\le 1.\] Then the map
\[H^0(\alpha(\Cot_{(R/p)/\setF_p})[-p]\otimes ^L k):H^0(\bigwedge^p\Cot_{(R/p)/\setF_p}[-p]\otimes^L_{R/p}k)\to H^0((F^*_{R/p}\Cot_{(R/p)/\setF_p}[-1])\otimes^L_{R/p}k)\] is an isomorphism by the above identification. This and the injectivity of $H^0(\kappa_R[-1]\otimes^L k)$ imply the injectivity of $H^0(c_R[-p]\otimes^Lk)$. As recalled in the first paragraph, this is equivalent to the assertion (2).

$(2)\Rightarrow (1)$: Conversely, assume (2), which is equivalent to the injectivity of the map $H^0(c_R[-p]\otimes^Lk)$. In particular, we have
\[\dim_k\Gamma^p(H_1(\Cot_{(R/p)/\setF_p}\otimes^L_{R/p}k))\le 1,\] which implies 
\[\dim _kH_1(\Cot_{(R/p)/\setF_p}\otimes^L_{R/p}k)\le 1.\]
By the commutativity of the above diagram of fiber sequences, this implies that the map $H^0(\kappa_R[-1]\otimes^Lk)$ is also injective. This, together with the lower fiber sequence in the above commutative diagram, implies the vanishing of $H_1(\FWCot_{R}\otimes_{R/p}^Lk)$. Since the vanishing of homology groups of $\FWCot_{R}\otimes_{R/p}^Lk$ in degree $\ge 2$ always holds by our assumption on $R$, we get the assertion (1).
\end{proof}
As expected by Bhatt and Mathew in the paragraph right before \cite[Lemma 4.18]{bhatt2023Syntomic}, there is a relationship between the regularity criterion via $\Cofib(\theta_R)$ and the regularity criterion via the Frobenius--Witt differentials:
\begin{remark}
Let $(R,\mathfrak m,k)$ be a $p$-torsionfree complete intersection local noetherian ring with residue characteristic $p$.
If $\FWCot_R$ is flat over $R/p$, then \Cref{theorem2.10} implies that $\Cofib(\theta_R)\otimes_{R/p}^Lk$ is concentrated in cohomological degree $\ge 1$. Therefore, \cite[Proposition 4.19]{bhatt2023Syntomic} implies the regularity of $R$. Thus, under the assumption that $R$ is a complete intersection, \cite[Proposition 4.19]{bhatt2023Syntomic} implies the regularity criterion via the Frobenius--Witt cotangent complex (i.e. the former part of \cite[Theorem 5.2]{mao2026FrobeniusWitt}).
\end{remark}

\end{document}